\documentclass[11 pt, psamsfonts]{amsart}
\usepackage{amssymb}

\def\ignore#1{\relax}

\def\g{\mathfrak g}

\def\choose #1 #2{\begin{pmatrix}#1\\#2\end{pmatrix}}
\def\R{{\mathbb R}}
\def\Z{{\mathbb Z}}
\def\Q{{\mathbb Q}}
\def\C{{\mathbb C}}
\def\la{\lambda}
\def\al{\alpha}
\def\La{\Lambda}

\def\S{\mathcal S}

\def\N{\mathbb N}

\def\M{\mathcal M}
\def\B{{\mathcal B}}

\def\om{\omega}
\def\lan{\langle}
\def\ra{\rangle}
\def\half{\frac{1}{2}}
\def\ep{\epsilon}

\def\Vn{V^{\otimes n}}

\def\rar{\rightarrow}
\calclayout

\def\th@plain{%
  \let\thmhead\thmhead@plain \let\swappedhead\swappedhead@plain
  \thm@preskip.5\baselineskip\@plus.2\baselineskip
                                    \@minus.2\baselineskip
  \thm@postskip\thm@preskip
  \itshape
\renewcommand{\labelenumi}{{(\alph{enumi})\quad}}
                        \renewcommand{\labelenumii}{{(\roman{enumii})\ }}
}
\def\th@definition{%
  \let\thmhead\thmhead@plain \let\swappedhead\swappedhead@plain
  \thm@preskip.5\baselineskip\@plus.2\baselineskip
                                    \@minus.2\baselineskip
  \thm@postskip\thm@preskip
  \upshape
}
\def\th@remark{%
  \thm@headfont{\itshape}% heading font bold
  \let\thmhead\thmhead@plain \let\swappedhead\swappedhead@plain
  \thm@preskip.5\baselineskip\@plus.2\baselineskip
                                    \@minus.2\baselineskip
  \thm@postskip\thm@preskip
  \upshape
}

{\theoremstyle{definition}
\newtheorem{theorem}{Theorem}[section]
}

{\theoremstyle{plain}
\newtheorem{proposition}[theorem]{Proposition}
}
 {\theoremstyle{definition}

}

{\theoremstyle{definition}
\newtheorem{lemma}[theorem]{Lemma}
}

{\theoremstyle{plain}

}

{\theoremstyle{definition}
\newtheorem{definition}[theorem]{Definition}
}

{\theoremstyle{definition}
\newtheorem{example}[theorem]{Example}
}

{\theoremstyle{definition}
\newtheorem{remark}[theorem]{Remark}
}

{\theoremstyle{remark}

}

\renewcommand{\labelenumi}{{ \theenumi.}}
\renewcommand{\labelenumii}{{(\alph{enumii})}}

\begin{document}
%\frontmatter

\title{A note on tensor categories of Lie type $E_9$}

\thanks{\textit{Department of Mathematics, Indiana University, Bloomington, IN 47401 USA}\\
email: \texttt{errowell@indiana.edu}}
\author{Eric C. Rowell}

\maketitle
\begin{abstract}
We consider the problem of decomposing tensor powers of the fundamental
level 1 highest weight representation $V$ of the affine Kac-Moody
algebra $\g(E_9)$.  We describe an elementary algorithm for
determining the decomposition of the submodule of $\Vn$ whose
irreducible direct summands have highest weights which are maximal
with respect to the null-root. This decomposition is based on
Littelmann's path algorithm and conforms with the uniform
combinatorial behavior recently discovered by H. Wenzl for the
series $E_N$, $N\not=9$.
\end{abstract}
\section{Introduction}

While a description of the tensor product decompositions for
irreducible highest weight modules over affine algebras can be
found in the literature (see e.g. \cite{Gould}), effective
algorithms for computing explicit tensor product multiplicities
are scarce.  Some partial results in this direction have been
obtained by computing
characters (see e.g. \cite{King1}) and  by employing crystal bases (see
e.g. \cite{Shimo}) or the equivalent technique of Littelmann
paths. In this note we look at the particular case of the affine
Kac-Moody algebra associated to the Dynkin diagram $E_9$, with any
eye towards extending the results of \cite{Wenzl3}.

Let $V$ be the irreducible highest weight representation of the
$\g(E_N) \ N\geq 6$ with highest weight $\La_1$ corresponding to
the vertex in the Dynkin diagram furthest from the triple point.
For $N\not=9$ H. Wenzl \cite{Wenzl3} has found uniform
combinatorial behavior for decomposing a certain submodule
$\Vn_{new}$ of $\Vn$ using Littelmann paths \cite{Li}.
These submodules have the property that each irreducible
summand of $\Vn_{new}$ appears in $\Vn$ for the first time (for $N\leq
8$) or last time (for $N\geq 10$).  The degeneracy of the 
invariant form was an obstacle to including the
affine, $N=9$ case.

We extend Wenzl's combinatorial description to the case $N=9$ by
finding submodules $\M_n$ analogous to his $\Vn_{new}$.
Specifically, we look at the (full multiplicity) direct sum of
those submodules of $\Vn$ whose highest weights have maximal
null-root coefficient. Not surprisingly, these summands appear
\emph{only} in $\Vn$. The particular
utility of considering this submodule is that whereas decomposing
the full tensor power $\Vn$ into its simple constituents would
require an infinite path basis, only a finite sub-basis
(consisting of 200 straight paths) is needed to determine the
decomposition of $\M_n$.  Although this note was inspired by the
results of \cite{Wenzl3}, the module $\M_n$ appears so naturally
that this case may shed some light on the combinatorial behavior
described by Wenzl.

This paper is organized in the following way.  In Section \ref{defs} we
give the data and standard definitions for the Kac-Moody algebra
$\g(E_9)$.  Section \ref{paths} is dedicated to summarizing the general
technique of Littelmann paths, while in Section \ref{e9} we apply this
technique to the present case and present some new definitions.
Table \ref{notation} gives a glossary of notation for the
reader's convenience. All the lemmas we prove are
contained in Section \ref{lemmas}, and the main theorem and algorithm they lead to is
described and illustrated in Section \ref{main}.  We briefly mention a possible
application and a generalization in Section \ref{future}, as well as connections to Wenzl's
results.

I would like to thank H. Wenzl for bringing this problem to my
attention and for many useful discussions.
 \section{Notation and Definitions}\label{defs}
\begin{figure}\label{Dynkin}\begin{picture}(300,80)(50,30)
\put(118,50){\circle{4}} \put(120,50){\line(1,0){18}}
\put(140,50){\circle{4}} \put(142,50){\line(1,0){18}}
\put(162,50){\circle{4}} \put(164,50){\line(1,0){18}}
\put(184,50){\circle{4}} \put(186,50){\line(1,0){18}}
\put(206,50){\circle{4}} \put(208,50){\line(1,0){18}}
\put(228,50){\circle{4}} \put(228,52){\line(0,0){18}}
\put(228,72){\circle{4}} \put(230,50){\line(1,0){18}}
\put(250,50){\circle{4}} \put(252,50){\line(1,0){18}}
\put(272,50){\circle{4}} \put(115,35){1} \put(137,35){2}
\put(159,35){3} \put(181,35){4} \put(203,35){5} \put(225,35){6}
\put(247,35){7} \put(232,67){8} \put(269,35){0}
\end{picture}\caption{Dynkin Diagram of $E_9$}
\end{figure}
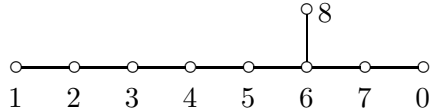

We begin by fixing a realization of the generalized Cartan matrix
of $\g(E_9)$ sometimes denoted in the literature by $\g(E^{(1)}_8)$.
Observe that our realization is different than that of Kac
\cite{Kac}. In particular the vertex Kac labels with a $0$ we
label with a $1$ in our Dynkin diagram (Figure \ref{Dynkin}). This
is done to conform with the notation of \cite{Wenzl3}.

\begin{definition}
Let $\{\ep_0,\delta,\ep_1,\ldots,\ep_8\}$
 be an ordered basis for $\R^{10}$, with symmetric bilinear form
$\lan\  ,\  \ra$ such that $\lan\ep_i,\ep_j\ra=\delta_{ij}$
 for $1\leq i,j\leq 8$
and $\lan\delta,\ep_0\ra=1$ with all other pairings $0$.
The simple
roots of $\mathfrak{g}(E_9)$ are defined by:
\begin{equation*}
\alpha_i=
\begin{cases} \ep_i-\ep_{i+1} & \text{if $1\leq i\leq 7$,}\\
\ep_7+\ep_8& \text{if $i=8$,}\\
\frac{1}{2}(\delta+\ep_8-\sum_{i=1}^7\ep_i)& \text{if $i=0$}
\end{cases}
\end{equation*}
\end{definition}
The simple roots generate the root lattice $Q:=\rm{span}_\Z\{\al_i\}$, and we define
coroots $\check\al_i:=\frac{2}{\lan \al_i,\al_i \ra}\al_i$.  As $\g(E_9)$
is simply-laced we abuse notation and identify each coroot with the
corresponding root.  Since we are only concerned with the combinatorics,
we refer the reader to Chapter 6 of the book \cite{Kac} for the full description of
the Kac-Moody algebra $\g(E_9)$.
\begin{definition}
We define the fundamental weights by:
\begin{equation*}
\La_i=
\begin{cases} (i,0;1,\ldots, 1,0,\ldots,0)\  \text{$i$ ones}
 & \text{if $1\leq i\leq 6$,}\\
\frac{1}{2}(8,0;1,\ldots,1,-1)& \text{if $i=7$,}\\
\frac{1}{2}(6,0;1,\ldots,1)& \text{if $i=8$,}\\
(2,0;0,\ldots,0)& \text{if $i=0$}
\end{cases}
\end{equation*}

Note that $\lan \alpha_i,\La_j\ra=\delta_{ij}$, and the
fundamental weights of $\mathfrak{g}(E_9)$ are determined up to a
multiple of $\delta$ by this relation. The set of dominant weights $P_+$ is
the $\N$-span of the fundamental weights plus $\C\delta$, and the $\Z$-span
$P$ is called the weight lattice.  It will be useful to denote by $\hat{P}_+$ those
dominant weights whose second coordinate is $0$.
\end{definition}
\begin{definition}
We define the Weyl group in the usual way: $W=\lan
s_i:i=0,\ldots,8\ra$ where $s_i(\nu)=\nu-\lan
\nu,\alpha_i\ra\alpha_i$ for $\nu \in \R^{10}$.
\end{definition}
The simple reflections $\{s_i\}_{i=1}^8$ generate a finite
subgroup $\overline{W}$ acting on the last eight coordinates by
permutations and an even number of sign changes. For any $\la$
with $\lan \la,\alpha_0\ra>0$ the vector $s_0(\la)$ has a strictly
smaller $\delta$-coefficient, and by applying elements of
$\overline{W}$ to arrange $\lan \la,\alpha_0\ra>0$ one can
construct an infinite sequence of vectors with strictly decreasing
$\delta$-coefficients.  Thus one sees that $W$ is an infinite
group.

\section{Littelmann Paths}\label{paths}
 To decompose the tensor powers of
$V$ we use the Littelmann path formalism (see \cite{Li}).
For this section we consider general Kac-Moody algebras
$\g$.

Littelmann considers the space of piecewise linear paths $\pi:
[0,1]\rightarrow P_\Q$ beginning at $0$ and ending at some point
in the weight lattice $P$.  He defines \emph{root operators} on
the space of all such paths $e_{i}$ and $f_{i}$ for each simple
root $\al_i$, which, when applied repeatedly to the straight path
$\pi_\la$ from $0$ to a dominant weight $\la$ give a \emph{path
basis} $\mathcal{B}_\la$ for the corresponding irreducible highest
weight module $V_\la$. The operators $f_i$ are defined on paths
$\pi$ as follows (see \cite{Li} for full details): let
$h_i(t)=\lan \pi(t),\al_i\ra$, and $m_i=\min(h_i(t))$.  If
$h_i(1)-m_i\geq 1$, split the interval $[0,1]$ into three pieces:
$[0,t_0]\cup [t_0,t_1]\cup [t_1,1]$ where $t_0$ is the maximal $t$
such that $h_i(t)=m_i$ and $t_1$ is the minimal $t$ such that
$h_i(t)=m_i+1$.  Then $$f_i\pi=\begin{cases}\pi(t) & \text{on
$[0,t_0]$}\\
s_i(\pi(t)) & \text{on $[t_0,t_1]$}\\
\pi(t)-\al_i & \text{on $[t_1,1]$}\end{cases}$$
If $h_i(1)-m_i<1$ then $f_i\pi=0$.  The operators
$e_i$ are defined similarly.  Since all paths begin at the weight
$0$, we may concatenate paths in the usual way.  For any $\la\in P$
define the path $\pi_\la: [0,1]\rightarrow P_\Q$ by
$t\rightarrow t\la$.  We denote concatenation by $*$, i.e.
$\pi_\la*\pi_\mu$ passes through $\la$ and terminates at $\la+\mu$.

Let $\la$ and $\mu$ be dominant weights of a Kac-Moody algebra
and $V_\la$, $V_\mu$ the corresponding irreducible highest weight
modules. We collect together those of Littelmann's results that we will need in:
\begin{proposition}
\begin{enumerate}
\item[(a)] $\B_\la\subset\{f_{1_j}f_{2_j}\cdots f_{s_j}\pi_\la\}$, that is, every
path in the basis $\B_\la$ is obtained from $\pi_\la$ by applying a finite sequence of
the root operators $f_i$.
\item[(b)]
The decomposition rules for the tensor product is given as
follows:
\[ V_\mu\otimes V_\la \cong \bigoplus_{\pi}V_{\pi(1)}
\] where $\pi=\pi_\mu *\pi_i$ with $\pi_i\in \mathcal{B}_\la$ and the image of $\pi$
contained in the closure of the dominant Weyl chamber.
\item[(c)] The multiplicity of $V_\nu$ in $\Vn_\la$ is
equal to the number of paths whose image is
contained in the closure of the
dominant Weyl chamber that terminated at $\nu$ that
are obtained by concatenating $n$ basis paths.
\end{enumerate}
\end{proposition}
\section{Lie type $E_9$}\label{e9}
Now we consider
the set $P(\La_1)$ of weights of $V$.  Following Kac \cite{Kac}
we call the weights in the Weyl group orbit $W\cdot\La_1$
\emph{maximal} and note that
$$P(\La_1)=\bigcup_{\om \in W\cdot\La_1}\{\om-t\delta : t \in \N\}.$$
Any $V_{\la-s\delta}$ that appears in some $\Vn$ must be of the
form:
$$\la-s\delta=\sum_{\om_i\in P(\La_1)}\om_i$$ with
$s\in\half\N$.  It is well known that the maximal weights appear in the 
multiset $P(\La_1)$ with
multiplicity one (for example, see \cite{Gould}).

 The second coordinate (essentially determined by the number of times $s_0$ occurs in a 
 minimal expression) 
 provides a gradation on
$W\cdot\La_1$ which motivates the following lemma, the proof of
which is a computation.
\begin{lemma}\label{forms} Every $\om\in W\cdot\La_1$ is of one of the
following 4 forms:
\begin{enumerate}
\item[(I)] Type I: $(1,0;\pm\ep_i)$
\item[(II)] Type II: $\frac{1}{2}(2,-1;\pm1,\ldots,\pm1)$
with an even number of minuses among the last eight coordinates.
\item[(III)] Type III: $(1,-1;w(1,1,1,0,\ldots,0))$
where $w \in S_8$, the group of permutations on $8$ symbols.
\item[(IV)] Type IV: all others, i.e. $(1,-j;\nu)$ where
$j\geq 1$ and if $j=1$, $\nu\not\in S_8\{(1,1,1,0,\ldots,0)\}$.
\end{enumerate}
\end{lemma}
The weights of types I, II and III will be particularly useful and
we will call them \emph{straight weights} and denote the set of
straight weights by $\Omega$. It is a simple but tedious
computation
 to show that for any $\omega\in\Omega$, the straight line path
 $\pi_\omega$ from $0$
 to $\omega$ is in the path basis $\mathcal{B}_{\La_1}$.
The idea of the computation is to start with the path $\pi_{\La_1}$
and inductively apply only those operators $f_i$ for which the height function
$h_i(t)=t$ so that two of the three intervals in the definition of the operator $f_i$
are degenerate, and the image of the paths remain straight lines.
Type I and II weights are in fact \emph{all} maximal weights with
second coordinate $0$ or $-\half$, while there are maximal weights
with second coordinate $-1$ besides those of type III.  Observe that since $\La_1$ is
the unique level one dominant weight (modulo $\delta$),
all paths obtained from concatenation of basis paths whose image lies in
the dominant Weyl chamber must pass through $\La_1$.
\begin{definition}\label{leveldef}
The \emph{level} $n(\la)$ of a weight $\la$ is the $\ep_0$
coordinate. Note that all weights in $P(\La_1)$ are level $1$,
thus $\la$ has level $n$ iff $\la-t\delta$ appears in $V^{\otimes
n}$ for some $t$ since $\la$ will be a sum of weights in
$P(\La_1)$.  Denote by $\hat{P}_+(n)$ the set of level $n$
dominant weights modulo $\delta$.
\end{definition}
The following definition appears in \cite{Wenzl3} and is critical
in the sequel.
\begin{definition}
We define the function $k: Q\rar\Z$ in one of the following
equivalent ways:
\begin{enumerate}
\item[(a)]  $k(\om)= -\lan\om,2\hat{\al}_0\ra$, where
 $\hat{\al}_0=\al_0-\ep_{8}$,
\item[(b)] If $\om=\sum_{i=0}^8 M_i\La_i$,
then $k(\om)=M_{8}-M_{7}-2M_0$.
\end{enumerate}
We will also need the quantity $[(\la)]_3$ defined to be the
remainder of $k(\la)$ upon division by $3$.
\end{definition}
We compute these values for the maximal weights and record them in
the following:
\begin{lemma}\label{omkvalues}
The values of the function $k$ for the maximal weights $\om$ of types I, II,
III and IV (as in Lemma \ref{forms}) satisfy:
\begin{enumerate}
\item[(I)] Type I: $k(\omega)\in\{0,-2\}$
\item[(II)] Type II: $k(\omega)\in\{3,1,-1,-3,-5\}$
\item[(III)] Type III: $k(\omega)=2$
\item[(IV)] Type IV: $k(\omega)\leq(6j-6)$ where $\om=(1,-j;\nu)$ and $1\leq j\in \half\Z$.
\end{enumerate}
\end{lemma}

The dominant weights are only defined up to a multiple
of $\delta$, but we are interested in those that appear in $P_+(\Vn)$ which motivates:
\begin{definition}
A level $n$ dominant weight $\la -m_\la\delta$ is called
\emph{initial} if $m_\la$ is minimal such that $V_{\la
-m_\la\delta}$ appears in $\Vn$.
\end{definition}
\begin{remark}  It is easy to see that there are
finitely many initial weights of a fixed level $n$, since there is
a one-to-one correspondence between the finite set $\hat{P}_+(n)$ and initial weights.
The term \emph{initial} comes from the fact that if $\la-m_\la\delta\in\hat{P}_+(n)$
so is $\la-(m_\la+1)\delta$.
Moreover, it is clear that $m_\la$ is always a non-negative
half-integer, since the coefficient of $\delta$ for any weight
$\om\in P(\La_1)$ is a non-positive half integer.
\end{remark}

We will eventually show that the $m_\la$ is computed from
the value of $k(\la)$ via the function:
\begin{definition}Let $\la\in\hat{P}_+$:
\begin{equation}
\Delta(\la)=
\begin{cases}
0 & \text{if $k(\la)\leq 0$ and even,}\\
\frac{1}{2} & \text{if $k(\la)\leq 1$ and odd,}\\
\frac{1}{6}(k(\la)+2[\la]_3)& \text{if $k(\la)\geq 1$}
\end{cases}
\end{equation}

Observe that when $k(\la)=1$ we have
$\frac{1}{6}(k(\la)+2[\la]_3)=\half$ so $\Delta$ is well defined.
\end{definition}
\begin{definition}
Define $\M_n$ to be the largest submodule of $V^{\otimes n}$ such
that all irreducible direct summands have highest weights of the form
$\la-m_\la\delta$ (i.e. initial weights).
\end{definition}
We illustrate this definition with an example:
\begin{example} The highest weight module $V_{\La_8}$ does not appear in
$V^{\otimes 3}$ as it is not a sum of $3$ type I weights. However,
$V_{\La_8-\frac{\delta}{2}}$ does appear in  $V^{\otimes 3}$ as
$$\La_8-\frac{\delta}{2}=\La_0+\frac{1}{2}(2,-1;1,\ldots,1)=
(1,0;\ep_1)+(1,0;-\ep_1)+\frac{1}{2}(2,-1;1,\ldots,1).$$  Notice
also that $V_{\La_8-\frac{t}{2}\delta}$ will also appear in
$V^{\otimes 3}$ for any $t\geq 1$, but only
$V_{\La_8-\frac{\delta}{2}}$ will appear in $\M_3$.\end{example}
The complete reducibility of $\Vn$ (see e.g. \cite{Gould}) allows
us to write:
$$\Vn\cong\M_n\oplus Z_n$$ where $Z_n$ consists of those simple submodules whose
highest weights are not initial.
\begin{table}\label{notation}
\caption{Notation}
\begin{tabular}{*{2}{lr}}
$\al_i$ & $i$th simple root\\
$Q$ & root lattice\\
$\La_i$ & $i$th fundamental weight\\
$P$ & weight lattice\\
$P_+$ & dominant weights \\
$\hat{P}_+$ & dominant weights $\pmod\delta$  \\
$\hat{P}_+(n)$ & level $n$ dominant weights $\pmod\delta$\\
$n(\la)$ & level of $\la$\\
$P(\La_1)$ & weights of $V$\\
$W\cdot\La_1$ & maximal weights of $V_{\La_1}$\\
$\Omega$ & set of straight weights\\
$P_+(\Vn)$ & dominant weights of $\Vn$\\
$[\la]_3$ & least residue of $k(\la)\pmod 3$\\
$\pi_\la$ &  path $t\rightarrow t\la$\\
$W$ & (affine) Weyl group\\
$\S(k)$ & level $k$ initial weights\\
$\la\rightarrow\mu$ & straight weight path\\
\end{tabular}
\end{table}
\section{Lemmas}\label{lemmas}
In this section we describe the combinatorial rules for
decomposing the modules $\M_n$.  The first two lemmas show that
$m_\la=\Delta(\la)$, while the two that follow show that one
may determine $\M_{n+1}$ from $\M_n$ and the (finitely many)
straight weights.
\begin{lemma}\label{exist}
Let $\la\in\hat{P}_+$ so that $\la-\Delta(\la)\delta$ is a level
$n$ dominant weight.  Then $V_{\la-\Delta(\la)\delta}$ appears in
$\Vn$. Moreover, there is a straight weight path from $0$ to $\la$
passing through only weights of the form $\mu -\Delta(\mu)\delta$
with $\mu\in\hat{P}_+$.
\end{lemma}
$Proof.$ Since we are not concerned with computing multiplicities,
it suffices to construct a piecewise linear straight weight path
from $0$ to $\la-\Delta(\la)\delta$ contained entirely within the
dominant Weyl chamber. Assume $\la=\sum_{i=0}^8 M_i\La_i$. We will
construct the required path in reverse by starting from the weight
$\la-\Delta(\la)\delta$ and removing path segments until we reach
the weight $0$. By concatenating the paths we remove in reverse
order we obtain the desired path.  The first set of useful paths
are the sub-paths of:
\[
\pi_1: 0\rar\La_1\rar\La_2\rar\La_3\rar\La_4\rar\La_5\rar\La_6\rar
(\La_7+\La_8)\rar(\La_0+2\La_8)
\] constructed by concatenating straight paths $\pi_\om$
terminating at $\om=(1,0;\ep_i)$, $i=1,\ldots,8$.  We denote by
$\pi_1^{(i)}$ the $i$th sub-path of $\pi_1$. In a similar fashion
we construct
\[ \pi_2:  0\rar\La_1\rar\La_0 \quad
\text{and} \quad \pi_3:0\rar\La_1\cdots\rar(\La_7+\La_8)\rar 2\La_7
\]
again using only paths with type I straight weight segments.  The
affect of removing these path segments on the value of $k(\la)$ is
as follows (where the value of a path at $1$ is $\la$):
\begin{enumerate}
\item $k(\pi*\pi_1^{(i)}(1))=k(\pi(1))$ i.e. deleting sub-paths of $\pi_1$ has
no affect on the value of $k(\la)$.
\item $k(\pi*\pi_2(1))=k(\pi(1))-2$ so deleting the path $\pi_2$ increases the value of
$k(\la)$ by 2.
\item $k(\pi*\pi_3(1))=k(\pi(1))-2$ so deleting the path $\pi_3$
increases the value of $k(\la)$ by 2.
\end{enumerate}
\subsection*{Case I} $k(\la)=M_8-M_7-2M_0\leq 0$ is even.\\
In this case $\Delta(\la)=0$.  Since $k(\la)$ does not
depend on $M_i$ $1\leq i\leq 6$ we can reduce the the case
$M_i=0$, $1\leq i\leq 6$ using sub-paths of $\pi_1$.  For $\la$
with $M_8=M_7=M_0=0$ we are done.  If not, we observe that $M_7$
and $M_8$ have the same parity. Again using the sub-path of
$\pi_1$ terminating at $\La_7+\La_8$ as many times as is
necessary, we may assume either $M_7=0$ or $M_8=0$.

\subsubsection*{Case I.1} $M_7=0$.\\
In this case we have $M_8\leq 2M_0$ and $M_8$ even. So by removing
path segments $\pi_1$ as many times is as necessary we can reduce
to $M_8=0$ with $k(\la)=-2M_0$ unchanged.  At this point we are
left with the case $\la=M_0\La_0$, to which we remove the path
segments $\pi_2$ as many times as necessary to reduce to $0$.

\subsubsection*{Case I.2} $M_8=0$.\\
Here we have that $M_7\geq -2M_0$ and $M_7$ is even, so we reduce by $\pi_3$
until $M_7=0$ and then reduce by the path $\pi_2$ until $M_0=0$ and we are left with
the weight $0$.  Observe that $k(\la)=-M_7-2M_0$ in this case so while deleting path segments
$\pi_2$ or $\pi_3$ result in a raised $k$-value,
it will always be non-positive and even, regardless.

\subsection*{Case II} $k(\la)=M_8-M_7-2M_0\leq 1$ is odd.\\
Here $M_8$ and $M_7$ have opposite parity; so, as in Case I, we
reduce by sub-paths of $\pi_1$ until either $M_7=0$ or $M_8=0$.  Then we reduce by
paths as in Case I
until we are left with two cases: $\la=\La_7-\frac{\delta}{2}$ and
$\la=\La_8-\frac{\delta}{2}$.  These are achieved by the paths:
$$0\rar\La_1\rar\La_2\rar\La_3\rar (\La_7-\frac{\delta}{2}) \quad
\text{and} \quad 0\rar\La_2\rar (\La_8-\frac{\delta}{2})$$ using
the straight weights $\half(2,-1;-1,-1,-1,1,\ldots,1,-1)$ and
$\half(2,-1;-1,-1,1,\ldots,1)$ respectively.

\subsection*{Case III} $k(\la)=M_8-M_7-2M_0\geq 2$.\\
In this case $M_8\geq 2+M_7+2M_0$, so we can use sub-paths of
$\pi_1$ to reduce to $M_7=0$ and then $M_0=0$ without changing the
value of $k(\la)$ and we are left with the task of constructing a
path terminating at $\la=M_8\La_8-\Delta(M_8\La_8)$, where
$k(\la)=M_8\geq 2$. The weight $3\La_8-\frac{\delta}{2}$ is of the
form $\mu-\Delta(\mu)\delta$ with $\mu\in\hat{P}_+$ and the path:
$$\pi_1*[(2\La_8+\La_0)\rar (3\La_8-\frac{\delta}{2})]$$
 allows us to reduce to $M_8\leq 2$ since:
\begin{equation}\label{three}
\Delta((M_8-3\ell)\La_8)=\frac{1}{6}(M_8-3\ell+2[(M_8-3\ell)\La_8]_3)=
\Delta(M_8\La_8)-\frac{\ell}{2}
\end{equation}
 so that
$$M_8\La_8-\Delta(M_8\La_8)\delta-\ell(3\La_8-\frac{\delta}{2}))=
(M_8-3\ell)\La_8-\Delta((M_8-3\ell)\La_8)\delta.$$  Now the cases
$M_8=0,1$ were covered in cases I and II respectively, so we need
only construct a path to $2\La_8-\delta$.  But this is nothing
more than a doubling of the path $$0\rar\La_2\rar
(\La_8-\frac{\delta}{2})$$ constructed above. This completes the
proof.   $\Box$

\begin{lemma}\label{minimal}
$\Delta(\la)=m_\la$ for all $\la\in\hat{P}_+$.
\end{lemma}
$Proof.$  By Lemma \ref{exist} it is sufficient to show that
$m_\la\geq\Delta(\la)$ since $\la-\Delta(\la)\delta$ appears in $\Vn$ hence
$m_\la\leq\Delta(\la)$.
Again we consider cases.
\subsection*{Case I} $k(\la)\leq 0$ and even.\\
Since $m_\la\geq 0=\Delta(\la)$ there is nothing to prove.
\subsection*{Case II} $k(\la)\leq 1$ and odd.\\
We need only show that $m_\la\not=0$.  The only way that this can
occur is if $\la$ can be expressed as a sum of type I weights. But
$k(\om)=0$ or $-2$ for $\om$ of type I, so if $\la$ were a sum of
type I weights $k(\la)$ would be even.
\subsection*{Case III} $k(\la)\geq 2$.\\  In this case we will
reduce to the case where $\la=M_8\La_8-t\delta$ using sub-paths of
$\pi_1$ as in the proof of Lemma \ref{exist}.  Suppose
$$\la=\sum_{i=0}^8M_i\La_i-t\delta$$ with $t$ minimal and
$t\in\half\N$. We compute $k(\la)=M_8-M_7-2M_0\geq 2$, so that
$M_8>M_7+2M_0$.  We can reduce to the case where
$M_1=M_2=\cdots=M_6=M_7=0$ using sub-paths of $\pi_1$ and
observing that
$$\la^\prime=M_0\La_0+(M_8-M_7)\La_8-t\delta$$ has
$k(\la^\prime)=k(\la)$ and $t$ minimal for $\la^\prime$ if and
only if $t$ is minimal for $\la$.  Setting $M_0^\prime=M_0$ and
$M_8^\prime=(M_8-M_7)$ we have
$k(\la^\prime)=M_8^\prime-2M_0^\prime\geq 2$.  Reducing by the
path $\pi_1$ and setting
$M_8^{\prime\prime}=M_8^\prime-2M_0^\prime$ we see that
$$\la^{\prime\prime}=M_8^{\prime\prime}\La_8-t\delta$$
has $t$ minimal if and only if $t$ is minimal for $\la^\prime$. So
we are left with showing that $m_\la\geq\Delta(\la)$ for
$\la=M_8\La_8$.  This will follow by an induction argument once we
show it for the cases $M_8=1,2$ and $3$.
\subsubsection*{$M_8=1$} This case was already covered in Case II
above.
\subsubsection*{$M_8=2$} If $2\La_8$ were a sum of type I
weights we would have $k(2\La_8)\leq 0$ so we must have at least
one weight of type II, III or IV. By considering the values of $k$
on these weights we see that that $t=1$ is minimal.
\subsubsection*{$M_8=3$} Again considering the values of $k$ we
see that type I weights are not sufficient and that $t=\half$ is
minimal.

Observing that
$\Delta((M_8-3\ell)\La_8)+\frac{\ell}{2}=\Delta(M_8\La_8)$ (see
Equation \ref{three}) the case $M_8>3$ follows by induction and we
are done.
   $\Box$

\begin{remark}
We may now redefine \emph{initial weight} to be any dominant
weight of the form $\la-\Delta(\la)\delta$, and we denote the set
of initial weights of level $n$ by
$\S(n)=\{\la-\Delta(\la)\delta:\la\in\hat{P}_+(n)\}$.
\end{remark}
\begin{lemma}If $\la-\Delta(\la)\delta$ is an initial weight
and $\om\in\Omega$ then $\la-\Delta(\la)\delta -\om$
is either an initial weight or not in the dominant Weyl chamber.
\end{lemma}
$Proof.$  Let $\mu-t\delta=\la-\Delta(\la)\delta-\om$ where
$\mu\in\hat{P}_+$.  Assume that $\mu$ is
in the dominant Weyl chamber. We must demonstrate that $t=\Delta(\mu)$. By
Lemma \ref{minimal}, we have that $t\geq\Delta(\mu)$ as $t\leq\Delta(\mu)$
would contradict the minimality of $\Delta(\mu)$. Using Lemma \ref{omkvalues}
we have:
\begin{equation}
t=
\begin{cases}
\Delta(\la) & \text{if $\om$ is type I,}\\
\Delta(\la)-\half & \text{if $\om$ is of type II,}\\
\Delta(\la)-1 & \text{if $\om$ is of type III}
\end{cases}
\end{equation}
It is sufficient to show that $\Delta(\mu)\geq t$ for all of these
cases.  We organize them by considering the value of $k(\la)$ as
follows:
\subsection*{Case I} $k(\la)\leq 0$ and even.\\
Here we have $\Delta(\la)=0$.  Since $\Delta(\mu)\geq 0$ the only
possibility is that $\om$ is of type I, for which it is clear.
\subsection*{Case II} $k(\la)\leq 1$ and odd.\\
The only possibilities are $\om$ of type I or II, since
$\Delta(\la)=\half$ in this case and $\Delta(\mu)\geq 0$.  If
$\om$ is of type II, then $t=0$, hence $\Delta(\mu)\geq t$ is
obvious.  If $\om$ is of type I, then $t=\half$ and Lemma
\ref{omkvalues} implies $k(\mu)=k(\la)-k(\om)\leq 3$ and odd.  If
$k(\mu)\leq 1$ and odd then $\Delta(\mu)=\half$ and we are done.
Otherwise $k(\mu)=3$, and we compute $\Delta(\mu)=\half$ as
required.
\subsection*{Case III} $k(\la)\geq 2$.\\
Here there are 3 cases depending on the type of $\om$.  The
computations are somewhat tedious, but straightforward.
\subsubsection*{Case III.1} $\om$ is of type I.\\
If $k(\om)=0$ then $k(\mu)=k(\la)$ hence $\Delta(\mu)=\Delta(\la)$
and we are done.  If $k(\om)=-2$, then $k(\mu)=k(\la)+2$ and we
must check the three 3 cases corresponding to the values of
$[\mu]_3$ (depending on $[\la]_3$) by evaluating
$\Delta(\mu)=\frac{1}{6}(k(\mu)+2[\mu]_3)$.

\subsubsection*{Case III.2} $\om$ is of type II.\\
We must show that $\Delta(\mu)\geq \Delta(\la)-\half$.  This is
the most involved case as $k(\om)\in\{3,1,-1,-3,-5\}$ and we must
check a total of 15 subcases corresponding to the 3 values of
$[\la]_3$ and 5 values of $k(\om)$. As an example of what is
involved we work out the cases where $[\la]_3=1$ and $k(\om)=-1$.
Then $k(\mu)=k(\la)+1$, $[\mu]_3=2$ and
$$\Delta(\mu)=\frac{1}{6}(k(\la)+1+2[\mu]_3)=
\frac{1}{6}(k(\la)+1+4)\geq\frac{1}{6}(k(\la)+2)=\Delta(\la)>\Delta(\la)-\half.$$
Notice that in this case $\mu$ is not dominant. The remaining
cases are handled similarly.
\subsubsection*{Case III.3} $\om$ is of type III.\\
We must show that $\Delta(\mu)\geq \Delta(\la)-1$.  Here
$k(\om)=2$ so $k(\mu)=k(\la)-2$ and we must again check cases by
evaluating $\Delta(\mu)$.   $\Box$

\begin{lemma}If $\la-\Delta(\la)\delta$ is an initial weight and
$\la-\Delta(\la)\delta +\om$ is also initial, then $\om$ is a
straight weight.
\end{lemma}
Before giving a proof, we mention a \emph{caveat}: the requirement
that $\la-\Delta(\la)\delta +\om$ is initial is not superfluous.
For example, $2\La_8-\delta$ is an initial weight and
$(1,0;-\ep_8)$ is a straight weight, but
$2\La_8-\delta+(1,0;-\ep_8)=\La_7+\La_8-\delta$ is not initial.\\
$Proof.$ It is enough to show that $\la-\Delta(\la)\delta+\om$ is
not initial if $\om$ is not straight.  The key fact here is from
Lemma \ref{omkvalues}: $k(\om)\leq 6j-6$ for the type IV weight
$\om=(1,-j;\nu)$ where $j\geq 1$ is a half-integer. Let
$\mu-s\delta=\la-\Delta(\la)\delta+\om$ for such a weight $\om$.
Observing that $s=\Delta(\la)+j$ we will show that
$s\not=\Delta(\mu)$.
\subsection*{Case I} $k(\mu)\leq 1$.\\
Since $s\geq j\geq 1$ and $\Delta(\mu)\leq \half$, it is clear
that $s\geq\Delta(\mu)$.
\subsection*{Case II} $k(\la)\leq 1$.\\
Here we have that
$$\Delta(\mu)=\frac{1}{6}(k(\la)+k(\om)+2[k(\la)+k(\om)]_3)\leq
\frac{1}{6}(1+6j-6+4)=\frac{6j-1}{6}<j\leq s$$ so once again
$\Delta(\mu)\not=s$.
\subsection*{Case III} $k(\la)\geq 1$ and $k(\mu)\geq 1$.\\
Computing as above we have
$$\Delta(\mu)=\frac{1}{6}(k(\la)+k(\om)+2[k(\la)+k(\om)]_3)\leq
\frac{1}{6}(k(\la)+6j-6+4)\leq\frac{1}{6}(k(\la))+\frac{6j-2}{6}<
\Delta(\la)+j=s$$
So we see that $\Delta(\mu)\not=s$ in all cases and we are done.
  $\Box$

\section{The Main Theorem and an Algorithm}\label{main}
The following theorem is a immediate corollary of the
lemmas in the previous section:
\begin{theorem}\label{summarized}
If $\la-\Delta(\la)\delta\in\S(n)$ then
any straight weight path from $0$ to $\la-\Delta(\la)\delta$
passes through only initial weights.  Thus
$$\M_n\cong\bigoplus_{\la\in\hat{P}_+(n)}c_\la V_{\la-\Delta(\la)\delta}$$
where the multiplicities $c_\la$ are determined by counting the
straight weight paths terminating at $\la-\Delta(\la)\delta$.
\end{theorem}
Applying the results we have the following simple inductive
algorithm for decomposing $\M_n$ as a sum of simple highest weight
modules:
\begin{enumerate}
\item[{\bf Step 1}] Initialize with $\M_1\cong V_{\La_1}$.
\item[{\bf Step 2}] Having determined the multiplicities $c_\la$
so that
$$\M_n\cong\bigoplus_{\la\in\hat{P}_+(n)}c_\la V_{\la-\Delta(\la)\delta}$$ compute the
set $A_\la=\{\la-\Delta(\la)\delta +\om: \om \in\Omega\}$ for each
$\la\in\hat{P}_+(n)$.
\item[{\bf Step 3}] Compute the set $\S(n+1)$.  The size of $\S(k)$
is computed from the generating
function:
\begin{equation}\label{genfun}
\prod_{0\leq i\leq 8}\frac{1}{1-x^{n(\La_i)}}=
1+x+3x^2+5x^3+10x^4+15x^5+27x^6+39x^7+63x^8+O[x^9]\end{equation}
\item[{\bf Step 4}] For each $\mu-\Delta(\mu)\delta\in\S(n+1)$, let
$B_\mu=\{\la\in\hat{P}_+(n): \mu-\Delta(\mu)\delta\in
A_\la\}$. Then
$$c_\mu=\sum_{\la\in B_\mu}c_\la.$$
\begin{remark} The formula in Step 3 is valid since
the level of a dominant weight $\la$ is determined by the
decomposition $\la=\sum_i M_i\La_i$ and the levels $n(\La_i)$ of
the fundamental weights $\La_i$ (see Definition \ref{leveldef}).  
One identifies a level $n$
dominant weight with a partition of $n$ into parts whose sizes are
in the multi-set $\{n(\La_i)\}$, and standard combinatorics lead to Eq. \ref{genfun}.  
For arbitrary $N$ the highest weight
module $V_\la$ appears in $V^{\otimes n(\la)}$ where the formula
for $n(\la)$ is given in \cite{Wenzl3}, Eq. (3.1) in case $N\not=9$.  
However, his formula  
breaks into three cases which depend on $k(\la)$ in a way that makes the
problem of constructing a generating function valid for all $N$ rather 
complicated combinatorially.
\end{remark}
\end{enumerate}
As an application we compute the decompositions of the first few
$\M_n$:
\begin{eqnarray*}
 &\M_2\cong&V_{\La_0}\oplus V_{\La_2}\oplus
 V_{2\La_1}\\
&\M_3\cong&V_{3\La_1}\oplus 2 V_{\La_1+\La_2}\oplus 3
V_{\La_0+\La_1}\oplus V_{\La_3}\oplus 2
V_{\La_8-\frac{\delta}{2}}\\
&\M_4\cong&V_{4\La_1}\oplus V_{\La_4}\oplus
6V_{\La_0+2\La_1}\oplus 3V_{\La_2+2\La_1}\oplus
6V_{\La_7-\frac{\delta}{2}}\oplus \\
& &6V_{\La_0+\La_2}\oplus 3V_{2\La_0}\oplus
8V_{\La_1+\La_8-\frac{\delta}{2}}\oplus 3V_{\La_1+\La_3}\oplus
2V_{2\La_2}
\end{eqnarray*}

\section{Connections and Further Directions}\label{future}
\subsection{$E_N$ Series}
Wenzl introduces a generic labeling set, $\Gamma$, for
the dominant integral weights of $\g(E_N)$, $N\not=9$ consisting
of triples $(n,\mu,i)$ where $n\in\N$, $\mu$ a Young diagram with
$|\mu|\leq n$ and $i\in\{0,1,2\}$, subject to some further
conditions (see \cite{Wenzl3}, Section 2).  The labeling is
realized via a map $\Phi$ assigning an element of $\Gamma$ to each
integral dominant weight.  The ambiguity in the dominant weights
due to the null-root precludes extending $\Phi$ directly to the
excluded case; however, the set of integral dominant weights of
$\g(E_9)$ whose image under $\Phi$ is in $\Gamma$ is precisely the
set of initial weights!  Thus one sees that our submodule $\M_n$
must be the ``missing link'' replacing $V_{new}^{\otimes n}$
required to extend Wenzl's main combinatorial result
(\cite{Wenzl3}, Proposition 3.10) for the $E_N$ $N\geq 6$ series
to the $N=9$ case:
\begin{proposition}
Assume $N>n$.  Then the branching rules for $V_{new}\subset
V^{\otimes 2}_{new}\subset\cdots\subset V^{\otimes n}_{new}$ do
not depend on $N$.
\end{proposition}
This proposition implies that when $k<9$ and $N>k$ the
combinatorial formula given is Step 3 of the algorithm holds.

\subsection{Braid Representations}
For generic $q$, the tensor product rules for the quantum group
$U_q\g(E_N)$ are the same as those of the Kac-Moody algebra
$\g(E_N)$. Wenzl was also able to show that, for $N\not=9$, the
centralizer algebra of the corresponding $U_q\g(E_N)$-module
$\Vn_{new}$ is generated by the image of the braid group $B_n$
(acting by $R$-matrices) and one more operator called the
\emph{quasi-Pfaffian}.  It should be possible to extend this
result to the $N=9$ case using the quantum group version of the
modules $\M_n$ together with the specific knowledge of the
decomposition rules.

\subsection{Other Lie Types}
It may be possible to use the same approach to derive a similar
algorithm for decomposing the tensor powers of low-level highest
weight modules for any affine Kac-Moody algebra.  By defining the
submodules analogous to $\M_n$ one would just need to determine
the subset of maximal weights corresponding to the set $\Omega$ of
straight weights.

\end{document}